\documentclass[12pt,a4paper,twoside]{scrartcl}
\usepackage[left=25mm,right=25mm,top=30mm,bottom=35mm]{geometry}
\usepackage[breaklinks=true,colorlinks=true,linkcolor=blue,urlcolor=black,citecolor=blue]{hyperref}%Hyperlinks
\usepackage[utf8]{inputenc}%umlauts
\usepackage[T1]{fontenc}%westeuropean coding
\usepackage{lmodern}%font
\usepackage[english]{babel}%language
\usepackage{indentfirst}%indent first
\usepackage{enumitem}%margins for itemize
\usepackage{csquotes}%for Quotes
\usepackage{endnotes}
\usepackage{ragged2e}
\let\footnote=\endnote
\begin{document}
\pagestyle{empty}
\begin{center}
	\hrulefill\newline
	\begin{Large}
		\textsc{The Author of a Quotation Goethe\\ Adduced Against Newton}\\
	\end{Large}
	\hrulefill
\end{center}

\begin{center}
	\begin{large}
		\textbf{Hubert Kalf}
	\end{large}\\[1em]
	Mathematisches Institut der Ludwig-Maximilians-Universität München,\\ 
	Theresienstraße 39, 80333 Munich, Germany\\[1em]
	Email: \texttt{kalf@mathematik.uni-muenchen.de}
\end{center}

\hspace*{1cm}

Goethe's lifelong antagonism to Newton is very well known and the literature about his Theory of Colour (beginning with [\hyperlink{CE}{CE}], this -- or its plural form -- has become the standard English translation of ``Farbenlehre'') is vast and overwhelming. For a brief and precise account of the differences between Goethe's chromatics, laid down in the first, the Didactic Part of his Farbenlehre of 1810, and Newton's mathematical optics we refer the reader to the Dutch oculist Halbertsma \cite[pp.~60--72]{10}.

The second part, called polemic, is an unpleasant personal dispute with Newton. The final part, translating as ``Materials for a History of Colour'', presents a history of the subject, naturally biased and sometimes also polemic, from the Graeco-Roman beginnings until Goethe's own time. As far as we know, these parts were never translated into English. (In the sequel all translations will be mine unless stated otherwise.)

In the Preface to his Didactic Part Goethe discloses the origin of his ambivalence, if not aversion, to a mathematical approach to science when he writes:
\begin{displayquote}
``$\dots$ thus it is possible to say that every attentive glance which we
cast on the world is an act of theorizing. This, however, ought to be done with consciousness, self-criticism, freedom, and, to use a daring word, with irony -- yes, all these faculties are necessary if abstraction, which we dread, is to be rendered innocuous, and the result which we hope for is to emerge with as much liveliness as possible.''\footnote{The translation given above is due to Erich Heller \cite[p.~19f.]{added}, but my attention was first drawn to the original by \cite[p.~311f.]{17}. The original reads [\hyperlink{LA}{LA I 4}, p.~5]: 
\begin{displayquote}
	$\dots$ so kann man sagen, daß wir schon bei jedem aufmerksamen Blick in die Welt theoretisieren. Dieses aber mit Bewußtsein, mit Selbsterkenntnis, mit Freiheit, und um uns eines gewagten Wortes zu bedienen, mit Ironie zu tun und vorzunehmen, eine solche Gewandheit ist nötig, wenn die Abstraktion, vor der wir uns fürchten, unschädlich und das Erfahrungsresultat, das wir erhoffen, recht lebendig und nützlich werden soll.	
\end{displayquote}
Unfortunately, the three English translations we are aware of mollify this sentence. [\hyperlink{CE}{CE}, p.~XXf.] and [\hyperlink{RM}{RM}, p.~72] write ``$\dots$ in order to guard against the possible abuse of this abstract view $\dots$'' and [\hyperlink{DM}{DM},  p.~159] has ``$\dots$ in order to avoid the pitfalls of abstraction $\dots$''.
}
\end{displayquote}
On the other hand, trying to explain to himself Newton's adherence to mathematics, Goethe surmises (in the section ``Newton's personality'' of his historical part)
\begin{displayquote}
``$\dots$ that Newton perhaps found so much pleasure in his theory because it offered new difficulties with each empirical step. In fact, a mathematician himself says: It is the habit of mathematicians to heap difficulties upon difficulties and incessantly find new ones in order to have the pleasure of surmounting them.''\footnote{\vspace*{-0.6cm}\begin{displayquote}
``\dots daß vielleicht Newton an seiner Theorie so viel Gefallen gefunden, weil sie ihm, bei jedem Erfahrungsschritte, neue Schwierigkeiten darbot. So sagt ein Mathematiker selber: C‘est la coutume des Géomètres de s‘éléver de difficultés en difficultés, et même de s‘en former sans cesse de nouvelles, pour avoir le plaisir de les surmonter.'' [\hyperlink{LA}{LA I 6}, p.~301]
\end{displayquote}
Note that at that time ``géomètres'' had a wider meaning than today.}
\end{displayquote}

Who was the mathematician who said this? As an expert on Euler, Andreas Speiser suggested in [\hyperlink{GA}{GA 17}, p.~920f.] that a similar idea can be found in Condorcet's eulogy on Euler \cite{4}. (The obituary in \cite{3} gives the reader an idea of Speiser's synoptic view of mathematics, art and philosophy.) For Condorcet (1743--1794), since 1776 permanent secretary of the French Academy of Sciences and maybe the last of the encyclopaedists, s.~\cite{9}. As far as I could make out the only passage that fits to some extent occurs at the end of his summary of Euler's achievements in partial differential equations. I am using (with some reservation) the English translation by John S.~D.~Glaus from the Euler Society, which is easily available from the net.
\begin{displayquote}
``The discovery of partial differentials, in reality belongs to Mr.~d’Alembert, since it is due to him that there is an understanding of the general form of their integrals. However, in Mr.~d’Alembert’s works we often saw the results of the calculations rather than the calculations themselves. It is to Mr.~Euler that we have the notation. He knew how to best present it, through the deep understanding of the theory which provided him the way in which to solve a great number of these equations, to distinguish the forms of the orders of integrals, for the different number of variables, to reduce these equations when they attain a certain form and become ordinary integrations, to provide for a way in which to remember these forms, through substitutions after which they vanish; in one word, to discover within the nature of these partial differential equations most of these singular properties which render the general theory so difficult and thorny, qualities which are nearly inseparable in Geometry where the degree of difficulty is so often the measurement of interest that one takes in a question and the honor that one attaches to a discovery. The influence of new truth in Science itself or to some application is the sole advantage that can balance the scale against the odds of such a difficult victory within the body of men for whom it is a pleasure to perceive a truth and one which is always proportional to the efforts that have been required.''\footnote{\vspace*{-0.6cm}\begin{displayquote}
C'est à \textsc{M.~D'Alembert} qu'appartient réellement la découverte du calcul aux différences partielles, puisque c'est à lui qu'est dûe la connaissance de la forme générale de leurs intégrales; mais dans les premiers Ouvrages de \textsc{M.~D'Alembert} on voyoit plus le résultat du calcul que le calcul lui-même; c'est à \textsc{M.~Euler} que l'on en doit la notation, il a su se le rendre propre, en quelque manière, par la profonde théorie qui l'a conduit à résoudre un grand nombre de ces équations, à distinguer les formes des intégrales pour les différens ordres et pour les différens nombres de variables, à réduire ces équations lorsqu'elles ont certaines formes, à des intégrations ordinaires, à donner les moyens de rappeler à ces formes, par d'heureuses substitutions, celles qui s'en éloignent; en un mot, en découvrant dans la nature des équations aux différences partielles, plusieurs de ces propriétés singulières qui en rendent la théorie générale si difficile et si piquante; qualités presque inséparables en Géométrie, où le degré de la difficulté est si souvent la mesure de l'intérêt qu'on prend à une question, de l'honneur qu'on attache à une découverte. L'influence d'une vérité nouvelle sur la Science même ou sur quelqu'application importante, est le seul avantage qui puisse balancer ce mérite de la difficulté vaincue, chez des hommes pour qui le plaisir d'apercevoir une vérité, est toujours proportionné aux efforts qu'elle leur a coûtés.
\end{displayquote}
\cite[p.~47]{4} and p.~294f.~of the reprint.}
\end{displayquote}

Euler is favourably mentioned in the historical part [\hyperlink{LA}{LA I 6}, p.~362] because, contradicting Newton, he asserted the possibility of achromatic lens systems, referring to the human eye \cite[p.~118]{5}. Euler's achromatics is also mentioned in the eloge \cite[p.~22f.]{4}, but there is nothing to indicate that Goethe consulted the eloge or anything else written by Condorcet for his History of Colour. (Cf.~\cite{14} where philologically oriented readers are primarily addressed.)

Goethe's quotation certainly stands on its own and to find its author one has to look elsewhere. I always thought that it had to be someone who critisised other mathematicians and was explicitly mentioned in the historical part, so that Goethe did not deem it necessary to cite him again. This leads to the historian of mathematics Jean-Etienne Montucla, who is twice somewhat curtly mentioned ([\hyperlink{LA}{LA I 6}, p.~239] in connection with Thomas Sprat, and on p.~352 in Goethe's ``Nachlese'' [gleanings]). Goethe attached particular importance to him, however, as can be seen from his entry
\begin{displayquote}
``To gain as much as can possibly be given to me when approaching mathematics, I read Montucla's Histoire $\dots$''\footnote{The original in Tag- und Jahreshefte 1806 [\hyperlink{LA}{LA II 6}, p.~329] reads: ``Um soviel als mir gegeben sein möchte, an die Mathematik heranzugehen, las ich Montuclas Histoire $\dots$''. 1806 was a particularly disastrous year. On October 14 Prussia was defeated in the battles of Jena and of Auerstedt, and Goethe was seriously endangered when French troops entered Weimar on the same day. The Duke of Brunswick, Gauss' patron, was fately wounded at Auerstedt.}
\end{displayquote}
He borrowed the two volumes from Jena in 1805, but did not start reading until 5 June 1806 (for more details s.~\cite{14}).\bigskip

Montucla's is the first history of mathematics, from the beginnings until around 1700, including its applications, which at the time were mechanics, astronomy, optics and mechanics \cite{16}. The first two volumes of a substantially enlarged 2$^{\mathrm{nd}}$ edition in four volumes appeared shortly before his death in 1799, the last two in 1801 -- primarily completed by the astronomer J.~J.~Lalande and by S.~F.~Lacroix, the author of an influential treatise on the calculus. A detailed description of the structure and a balanced critique of the two editions can be found in Sarton's richly documented paper \cite{18}. A photomechanical reprint, at places somewhat difficult to read, of the 2$^{\mathrm{nd}}$ edition appeared in 1960, giving Montucla the initials J.~F.~on its flyleaves (the first edition did not reveal its initials). Since Gallica provides access to \cite{16}, one can search for ``coutume'' (in T2; T1 is not relevant here). The 5$^{\mathrm{th}}$ of 6 hits gives the desired quotation verbatim (on p.~449 of the 1$^{\mathrm{st}}$ and p.~470 of the 2$^{\mathrm{nd}}$ ed.). Search via Google with ``Montucla'' and parts of the quotation turns out to be unsuccessful though.

The context into which our citation is somewhat surprisingly placed is the following.\footnote{Unconcerned with the question of authorship, Paul Epstein \cite[p.~88]{7} thought that the quotation aptly characterised ``the predominantly formalistic orientation'' in particular of the ``school of \underline{Lagrange} and \underline{Laplace} who served a cult of analytic formulae and wallowed in a sort of rage du calcul.'' This interpretation is at variance with Goethe's famous and frequently cited (by Epstein himself on p.~82) praise of Lagrange in number 609 of his ``Maxims and Reflections.'' An English translation can be found in [\hyperlink{DM}{DM}, p.~311].} In the Acta Eruditorum for May 1690 James Bernoulli (1654--1705) used the new infinitesimal calculus to solve a problem set by Leibniz (1646--1716) to challenge Cartesian dynamics and then solved by himself and by Huygens (1629--1695) geometrically. James ends his paper, possibly instigated by his brother John (1667--1748) \cite[p.~97]{2}, by initiating a new contest, viz., to determine the shape an inextensible (but flexible) homogeneous cord or chain will assume under its own weight when suspended freely from its ends. We follow Montucla's brief and frank outline of the history of this problem  as it gives us the opportunity to point out some inaccuracies that still persist in the literature. Starting on p.~446, Montucla writes: Galilei thought (``sans aucune raison solide'') that it was a parabola\footnote{With reference to a statement in the ``2$^{\mathrm{nd}}$ Day'' of Galilei's Discorsi this is the predominant view in today's literature as well (an exception is \cite{6}). It was not until 1972 that Fierz drew attention to the ``4$^{\mathrm{th}}$ Day'', where the parabolic trajectory is set into analogy with the catenary, whence it becomes clear that an approximation was meant (s.~\cite{11} and the literature cited there). The reason for this discrepancy seems to be unclear. Stillman Drake in his translation \cite[p.~146]{8} points to a lack of the usual conversational conclusion of the 2$^{\mathrm{nd}}$ Day and of a conversational opening of the 3$^{\mathrm{rd}}$ Day. So Galilei might have intended to add material here.} and that the Jesuit Fathers Pardies and Lana de Terzi had mistakenly tried to prove this.\footnote{The latter can certainly be disregarded \cite[p.~89]{21}. Pardies (1636--1673), however, in a book of 1673 proved -- in a slightly different way than Huygens (cf.~the next footnote) -- that the shape of the uniformly heavy cord is not a parabola \cite[p.~51, LXXV]{21}. The claim to the contrary in \cite[p.~229]{13} adduces a sentence by Pardies as well as a criticism of Pardies by James Bernoulli, but both concern different physical situations (see p.~52 [LXXVI] and p.~89 in \cite{21}). As an aside, Goethe mentions Pardies in his History of Colours [\hyperlink{LA}{LA I 6}, p.~267f.] because he led a small critical correspondence with Newton.} Joachim Jungius, however, established (``par diverses expériences'') that it was not a parabola, but he could not reveal the true nature of this curve.\footnote{Introductory to his solution in the Acta of 1691, Leibniz mentions such a result, but very likely it is irretrievably lost (Montucla's reference is erroneous). Huygens' celebrated arguments of 1646 \cite[p.~44ff.]{21}, which ruled out the parabola, were not published in his lifetime. Based on Huygens' idea, \cite{12} presents a simplification that could have been within the reach of Jungius (1587--1657).

As was common at the time, Jungius worked  in diverse areas, and Goethe wrote friendly of his botanical studies [\hyperlink{LA}{LA I 10}, pp.~285--296].} He then mentions that solutions were presented without proof (``apparently to leave some laurels to reap for those who managed to guess [the proof]''\footnote{``apparement afin de laisser encore quelque lauriers à cueillir, à ceux qui viendraient à bout de la deviner'' \cite[p.~447]{16}. Leibniz is more outspoken in a letter to his friend Rudolph Christian von Bodenhausen: ``In case one really exhibits something, it is good to give no proof at all or one where they do not find out our wiles.'' (For the German original and the reference s.~\cite[p.~71]{21}.)}) by the Bernoulli brothers and by Leibniz and Huygens in the Acta of 1691.\footnote{For the presentations Leibniz and John Bernoulli gave in this volume see \cite[pp.~69--72]{21} and \cite[pp.~72--75]{21}. Huygens' geometric solution was reconstructed for his collected papers by Korteweg \cite{15}. Truesdell confines himself to giving a difference approximation that leads to the decisive differential equation \cite[pp.~66ff.]{21}. Whether James Bernoulli himself had a proof in 1690 is unclear \cite[pp.~75,80]{21}.} Then he presents the proof John gave in 1691/92 in his private lectures to Marquis de l'Hospital when introducing him to the new calculus of Leibniz. (These had been published in Vol.~3 of John's Collected Papers \cite{1}.\footnote{The reader should not miss to have a look at the delightful frontispiece (the same for each volume). Below the inscription ``supra invidiam'' (beyond envy) a scroll dangles from a tree, displaying the construction of a curve (presumably the brachystochrone), and an upright dog tries to look closer.})

Listing several important properties of the catenary (notably the extremal property of its centre of gravity; cf.~\cite[p.~70]{21}), our quotation follows, since the Bernoullis immediately started to attack a number of more difficult problems. Both considered inextensible chains of varying thickness and extensible chains of fixed density.\footnote{This competition eventually degenerated into bitter animosity which was frequently and extensively discussed in the literature. We refer to the prudent analysis given in \cite{20}. On p.~262 Thiele cites a passage by Montucla quite near our quotation \cite[p.~458]{16}. (He works with the 2$^{\mathrm{nd}}$ edition, but the text is identical.)

As a curiosity we mention that Schillemeit (1931--2002) was able to clarify the authorship of another mathematical quotation with the help of Montucla's second edition \cite{19}. In 1954 Schillemeit wrote his thesis (Staatsexamensarbeit) on a problem, posed by Rellich, on the spectral theory of Dirac operators before changing to German philology with particular emphasis on Goethe and Kafka. As a rule such theses are never cited, but this one was \cite{22}.} Moreover, John determined the shape of a chain which hangs in a central gravitational field. (For a modern presentation of the solution of this problem by means of the calculus of variations we recommend the enjoyable paper \cite{6}.)

Thus the difficulties our quotation is referring to are difficulties that arise from ever more realistic descriptions of natural phenomena and not from artificial complications in a l'art pour l'art problem. Goethe saw this quite clearly (``new difficulties with each empirical step''), and he probably found it natural for topics which lend themselves to mathematics -- but not for his Farbenlehre, which he wanted to be seen as detached from optics [\hyperlink{LA}{LA I 3}, p.~229, 2 b 4].
%
%
%
%
%
%Endnotes:
\section*{Footnotes}
\vspace*{-2cm}
\theendnotes
%
%
%
%
%
%References
\section*{References}
\renewcommand\refname{\vskip -1cm}

The Leopoldina Edition of Goethe's scientific work appeared in 31 volumes from 1947 until 2019; it contains his Farbenlehre as well as all material pertaining to it. For example, LA I 6 is volume 6, containing the historical part, and LA II 6 refers to additional material or explanations. When possible, we have combined our page references with a cue, so that it should not be too difficult to find the passage in question in more readily available editions of Goethe's works such as 

\setlist[itemize]{leftmargin=1.5cm}
\begin{itemize}
	\item[{\hypertarget{GA}{[GA]}}]Gedenkausgabe der Werke, Briefe und Gespräche, 24 Bände und 3 Ergänzungen (Ernst Beutler Hrsg.). Artemis, Zürich 1949--72;
	\item[{\hypertarget{HA}{[HA]}}]Hamburger Ausgabe in 14 Bänden (Erich Trunz Hrsg.), Beck, München 1982.
\end{itemize}

Both omit the polemic part. A.~Speiser, in GA 16, collected the didactical and historical parts of the Farbenlehre with Goethe's writings on the Science of Knowledge (Wissenschaftslehre), adding an extensive introduction that is frequently discussed in the literature (e.g.~in \cite{17} below). HA 13, 14 contain these two parts in separate volumes with an introduction by C.~F.~von Weizsäcker in Vol.~13.\bigskip

It is only the didactic part that was translated into English and the reader has three choices.

\setlist[itemize]{leftmargin=1.5cm}
\begin{itemize}
	\item[{\hypertarget{CE}{[CE]}}]Goethe's Theory of Colours. Transl.~from the German with Notes by Charles Lock Eastlake. John Murray, London 1840. (Reprinted in very many editions which are still available.)
	\item[{\hypertarget{RM}{[RM]}}]Goethe's Color Theory. Arranged and edited by Rupprecht Matthei. Amer.~ed. translated and edited by Herb Aach. Van Nostrand Reinhold, New York 1971.
	\item[{\hypertarget{DM}{[DM]}}]Johann Wolfgang von Goethe, Scientific Studies, Edited and translated by Douglas Miller. Suhrkamp Publishers, New York 1988. (This is the last volume of Victor Lange's Selected Works of Goethe in 12 Volumes.)
\end{itemize}
%
%
%
%
%
%Bibliography:

%
%
%
%
%

\begin{thebibliography}{10}
\bibitem{1}
Bernoulli, Johann:~Opera Omnia t.~I--IV. G.~Gramer ed. Lausanne \& Genf 1742 (Reprint by Olms, Hildesheim 1968, with an introduction by J.~E.~Hofmann)

\bibitem{2}
Bernoulli, Johann:~Der Briefwechsel mit Jakob Bernoulli u.a. O.~Spiess ed. Basel 1955

\bibitem{3}
Burckhardt, J.~J.:~Andreas Speiser (1885--1970). In:~math.ch/100,129--161. Schweizerische Mathematische Gesellschaft 1910--2010 (Bruno Colbois et.~al.~eds.), GMS Press 2010

\bibitem{4}
Condorcet, Marquis de: Eloge de M.~Euler. Histoire de l’Academie royale des sciences pour l’anneé 1783, 37--68, Paris 1786. (Reprinted in:~L.~Euler, Opera omnia (3)12, 287--310. Orell Füssli, Zürich 1960.)

\bibitem{5}
Darrigol, O.:~ A History of Optics. University Press, Oxford 2012

\bibitem{6}
Denzler, J.~and Hinz, A.~M.:~Catenaria Vera -- the True Catenary. Expo.~Math 17, 117--142, 1999

\bibitem{7}
Epstein, P.:~Goethe und die Mathematik. Jahrbuch der Goethe-Gesellschaft 10, 76--102, 1924

\bibitem{8}
Galileo: Two New Sciences: Including Centers of Gravity and Force of Percussion. A New Translation. With Introduction and Notes by Stillman Drake. Univ.~of Wisconsin Press, Madison WI 1974.

\bibitem{9}
Granger, G.:~Condorcet, Marie-Jean-Antoine-Nicolas Caritat, Marquis de. In:~Dictionary of Scientific Biography. Vol.~3, 383--388 (C.~C.~Gillispie ed. Scribner, New York 1981)

\bibitem{10}
Halbertsma, K.~T.~A.:~A History of the Theory of Colour. Swets \& Zeitlinger, Amsterdam 1949

\bibitem{added}
Heller, E.:~The Disinherited Mind. Bowes \& Bowes. Cambridge 1952

\bibitem{11}
Herzig, A. and Szabó, I.:~Die Kettenlinie, das Pendel und die «Brachistochrone» bei \textsc{Galilei}. Verhandlungen der Naturforschenden Gesellschaft in Basel 91, 51--76, 1981 (also in the 3rd ed.~of I.~Szabó, Geschichte der mechanischen Prinzipien. Birkhäuser, Basel 1987)

\bibitem{12}
Hofmann, J.~E.:~Wie Jungius bewiesen haben könnte, daß die Kettenlinie keine Parabel ist. Sudhoffs Archiv 50, 302--305, 1966

\bibitem{13}
Hofmann, J.~E.:~Vom öffentlichen Bekanntwerden der Leibnizschen Infinitesimalmathematik, Sitzber.~Österr.~Akad.~Wiss.~Math.-Naturw.~Klasse 175, 209--254, 1966 (also in: Joseph Ehrenfried Hofmann.~Ausgewählte Schriften II, S.~1--46. Olms, Hildesheim 1990)

\bibitem{14}
Kalf, H.:~Der Urheber eines Mathematikzitates aus dem historischen Teil der Farbenlehre. Jahrbuch der Goethe-Gesellschaft. To appear.

\bibitem{15}
Korteweg, D.~J.:~La solution de Christiaan Huygens du problème de la chaînette. Bibliotheca Mathematica (3)1, 97--108, 1900

\bibitem{16}
Montucla, J.~E.:~Histoire des Mathématiques. T1--2. Paris 1758

\bibitem{17}
Neubauer, J.:~»Die Abstraktion, vor der wir uns fürchten«. Goethes Auffassung der Mathematik und das Goethebild in der Geschichte der Naturwissenschaft. In:~Versuche zu Goethe. Festschrift für Erich Heller. S.~305--320. Hrsg.~V.~Dürr, Heidelberg 1976

\bibitem{18}
Sarton, G.:~Montucla (1725--1799). His life and works. Osiris 1, 519--567, 1936

\bibitem{19}
Schillemeit, J.:~Der Geometer und die Dichtung. Philologische Arabeske über eine literarische Anekdote. In:~Studien zur Goethezeit, S.~541--560. Hrsg.~von Rosemarie Schillemeit. Wallstein, Göttingen 2006 (Erstdruck in der Festschrift für Victor Lange 1977)

\bibitem{20}
Thiele, R.:~Das Zerwürfnis Johann Bernoullis mit seinem Bruder Jakob. In:~Natur, Mathematik und Geschichte. Beiträge zur Alexander-von-Humboldt-Forschung und zur Mathematikhistoriographie, S.~257--276. H.~Beck u.a.~Hrsg.~Halle 1997

\bibitem{21}
Truesdell, C.:~The rational mechanics of flexible or elastic bodies. 1638--1788. In:~L.~Euler, Opera omnia (2)11.2. Orell Füssli, Zürich 1960

\bibitem{22}
Weidmann J.:~Oszillationsmethoden für Systeme gewöhnlicher Differentialgleichungen. Math.~Z. 119, 349--373, 1971
\end{thebibliography}
\end{document}